\documentclass[12pt,twoside]{article}
\usepackage{amsmath,amsbsy,amsfonts}
\usepackage{graphicx,color}
\newtheorem{theorem}{Theorem}[section]
\newtheorem{lemma}[theorem]{Lemma}
\newtheorem{proposition}[theorem]{Proposition}

\newtheorem{remark}[theorem]{Remark}
\newtheorem{claim}[theorem]{Claim}

\newcommand{\fim}{\hfill\rule{2mm}{2mm}}

\newcommand{\R}{\mathbb{R}}

\begin{document}

\setlength{\baselineskip}{4.5mm} \setlength{\oddsidemargin}{8mm}
\setlength{\topmargin}{-3mm}
\title{\Large\sf Infinitely many sign-changing solutions for
a class of elliptic problem with exponential critical growth.}
\author{\sf Denilson S. Pereira\thanks{denilsonsp@dme.ufcg.edu.br} \\
 Universidade Federal de Campina Grande\\
 Unidade Acad\^emica de Matem\'atica - UAMat\\
 CEP: 58.429-900 - Campina Grande - PB - Brazil\\ }

\pretolerance10000
\date{}
\numberwithin{equation}{section} \maketitle
\begin{abstract}
In this work we prove the existence of infinitely many nonradial solutions that change signal
to the problem $-\Delta u=f(u)$ in $B$ with $u=0$ on $\partial B$, where $B$ is the unit ball in $\R^2$ and $f$ is a continuous and odd function with exponential critical growth. \\


\vspace{0.3cm}


\end{abstract}

\section{Introduction}

Let $\Omega\subset\R^N$ be a bounded domain with smooth boundary and $f:\R\to\R$ be a $C^1$ function with $f(-t)=-f(t)$. Consider the following problem
$$
\left\{
\begin{array}{rcl} -\Delta u&=&f(u),~  \mbox{in}~~ \Omega, \\[0.2cm]
\mathcal{B}u&=&0,~\mbox{on}~~ \partial \Omega,
\end{array}
\right. \eqno (P)
$$
when $N\geq 4$, $\mathcal{B}u=u$ and $f(t)=|t|^{\frac{4}{N-2}}+\lambda t$, Br\'ezis-Niremberg \cite{BN} proved that $(P)$ admits a non-trivial positive solution, provided $0<f'(0)<\lambda_1(\Omega)$, where $\lambda_1(\Omega)$ is the first eigenvalue of $(-\Delta,H_0^1(\Omega))$. In \cite{CSS}, Cerami-Solimini-Struwe proved that if $N\geq 6$, problem $(P)$ admits a solution with changes sign. Using this, they also proved that when $n\geq 7$ and $\Omega$ is a ball, $(P)$ admits infinitely many radial solution which change sign.

Comte and Knaap \cite{CK} obtained infinitely many non-radial solutions that change sign
for $(P)$ on a ball with Neumann boundary condition $\mathcal{B}u=\frac{\partial u}{\partial\nu}$, for every $\lambda\in\R$. They obtained such solutions by cutting the unit ball into angular sectors. This approach was used by Cao-Han \cite{CH}, where the authors dealt with
the scalar problem $(P)$ involving lower-order perturbation and by de Morais Filho et al. \cite{MMF}  to obtain multiplicity results for a class of critical elliptic systems
related to the Br\'ezis-Nirenberg problem with the Neumann boundary condition on a ball.

When $N=2$, the notion of ``critical growth'' is not given by the Sobolev imbedding, but by the {\it Trudinger-Moser inequality} (see \cite{T} and \cite{M}), which claims that for any  $ u \in H^{1}_0(\Omega)$,
\begin{equation} \label{X0}
\int_{\Omega}
e^{\alpha u^2}dx
< +\infty, \,\,\,\, \mbox{ for every }\,\,\alpha >0.
\end{equation}
Moreover, there exists a positive constant $C=C(\alpha,|\Omega|)$ such that
\begin{equation} \label{X1}
\sup_{||u||_{H_0^{1}(\Omega)} \leq 1} \int_{\Omega} e^{\alpha u^2} dx \leq C , \,\,\,\,\,\,\, \forall \, \alpha  \leq 4 \pi .
\end{equation}

Motivated by inequality in $(\ref{X1})$, we say that the nonlinearity $f$ has exponential critical growth if $f$ behaves like $e^{\alpha_0s^2}$, as $|s|\to\infty$, for some $\alpha_0>0$. More precisely,
$$
\lim_{|s|\to\infty}\dfrac{|f(s)|}{e^{\alpha s^2}}=0,~~\forall\alpha>\alpha_0~~\mbox{and}~~\lim_{|s|\to\infty}\dfrac{|f(s)|}{e^{\alpha s^2}}=+\infty~~\forall\alpha<\alpha_0.
$$
In this case, Adimurthi \cite{Ad} proved that $(P)$ admits a positive solution, provided that $\displaystyle\lim_{t\to\infty}tf(t)e^{\alpha t^2}=+\infty$ (See also Figueiredo-Miyagaki-Ruf \cite{DOR} for a more weaker condition). In \cite{AY}, Adimurthi-Yadava proved that $(P)$ has a solution that changes sign, and when $\Omega$ is a ball in $\R^2$, $(P)$ has infinitely many radial solutions that change sign. Inspired in \cite{CK}, this paper is concerned with the existence of infinitely many non-radial sign changing solutions for $(P)$ when $f$ has exponential critical growth and $\Omega$ is a ball in $\R^2$. Our main result complements the studies made in $\cite{CK}$ and \cite{MMF}, because we consider the case where $f$ has critical exponential growth in $\R^2$. It is important to notice that in both studies mentioned above was
considered the Neumann boundary condition in order that the Pohozaev identity (see \cite{P}) ensures that the problem $(P)$ with the
Dirichlet boundary condition, has no solutions for $\lambda<0$ and $N\geq 3$. Since the Pohozaev identity is not available in dimension two, in our case we can use the Dirichlet boundary condition.

Here we suppose the following assumptions
\begin{enumerate}

\item[$(f_1)$] There is $C>0$ such that
$$
|f(s)|\leq Ce^{4\pi |s|^2}\ \ \mbox{for all}\ \ s\in\R;
$$

\item[$(f_2)$] $\displaystyle\lim_{s\rightarrow
0}\dfrac{f(s)}{s}=0$;

\item[$(H_1)$] There are $s_0>0$ and $M>0$ such that
$$
0< F(s):=\int_0^{s}f(t)dt\leq M|f(s)|,~ \mbox{for all}\ \|s|\geq s_0.
$$

\item[$(H_2)$] $0<F(s)\leq\dfrac{1}{2}f(s)s,\ \ \mbox{for all}\ \ s\in\R\setminus\{0\}.$

\item[$(H_3)$] $\displaystyle\lim_{s\to\infty} sf(s)e^{-4\pi s^2}=+\infty$
\end{enumerate}

Our main  result is the following:

\begin{theorem}\label{sb}
Let $f$ be an odd and continuous function satisfying $(f_1)-(f_2)$ and $(H_1)-(H_3)$. Then, problem $(P)$ has infinitely many sign-changing solutions.
\end{theorem}

\section{ Notation and auxiliary results}
\begin{figure}[!htb]
\begin{center}
\includegraphics[scale=0.3]{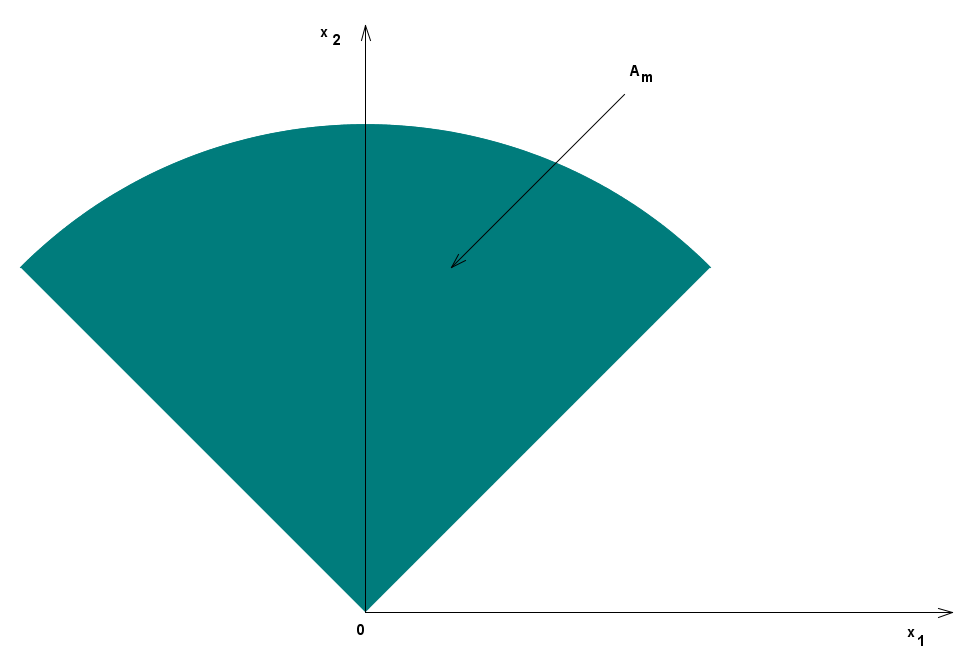}
\end{center}
\caption{Angular sector $A_m$.}
\label{dddd}
\end{figure}

For each $m\in\mathbb{N}$, we define
$$
A_m=\left\{x=(x_1,x_2)\in B:~~
cos\left(\dfrac{\pi}{2^m}\right)|x_1|<sen\left(\dfrac{\pi}{2^m}\right)x_2\right\}.
$$
So $A_1$ is a half-ball, $A_2$ an angular sector of angle $\pi/2$, and $A_3$ an angular sector of angle $\pi/4$ and so on (see figure \ref{dddd}).

Using the above notation, we consider the following auxiliary Dirichlet problem
$$
\left\{
\begin{array}{rcl} -\Delta u&=&f(u),~  \mbox{in}~~ A_m, \\[0.2cm]
u&=&0,~\mbox{on}~~ \partial A_m,
\end{array}
\right. \eqno (P)_m
$$
We will use the Mountain Pass Theorem to obtain a positive solution of $(P)_m$. Using this solution together with an anti-symmetric principle, we construct a sign-changing solution of problem $(P)$.

According to Figueiredo, Miyagaki and Ruf \cite{DOR}, to obtain a positive solution of $(P)_m$ it is sufficient to assume that the limit in $(H_3)$ verifies
$$
(H_3)'~~~\lim_{s\to+\infty}sf(s)e^{-4\pi s^2}\geq\beta>\dfrac{1}{2\pi d_m^2},
$$
where $d_m$ is the radius of the largest open ball contained in $A_m$. The hypothesis $ (H_3) $ was initially considered in Adimurthi \cite{Ad}. This hypothesis will be fundamental to ensure not only the existence but also the multiplicity of sign-changing solutions. As we will see bellow, assuming $(H_3)$ in place of $(H_3)'$, we have the existence of positive solution of $(P)_m$, for every $m\in\mathbb{N}$. This is the content of the next result.

\begin{theorem}\label{th1}
Under the assumptions $(f_1)-(f_2)$ and $(H_1)-(H_3)$, problem $(P)_m$ has a positive solution, for every $m\in\mathbb{N}$.
\end{theorem}

\section{Proof of Theorem \ref{th1}.}

Since we are interested in positive solutions to the problem $(P)_m$, we assume that
$$
f(s)=0,~~~\forall s\leq 0.
$$

Associated with problem $(P)_m$, we have the functional $I: H_0^1(A_m)\to\R$ given by
$$
I(u)=\dfrac{1}{2}\int_{A_m}|\nabla u|^2-\int_{A_m}F(u).
$$
In our case, $\partial A_m$ is not of class $C^1$. However, the functional
$I$ is well defined. In fact, for each $u\in H_0^1(A_m)$, let us consider
$u^{*}\in H_0^1(B)$ the zero extension of $u$ in $B$ defined by
$$
u^{*}(x)=\left\{
\begin{array}{cc}
u(x),& \mbox{if}~~~x\in A_m\\
0,& \mbox{if}~~~x\in B\setminus A_m.
\end{array}
\right.
$$
Clearly 
$$
\|u\|_{A_m}=\|u^{*}\|_{B}.
$$
Then, from $(f_1)$ and the Trudinger-Moser inequality (\ref{X0})
$$
\left|\int_{A_m}F(u)\right|=\left|\int_{B}F(u^{*})\right|\leq
\int_{B}|F(u^{*})|\leq C\int_{B}e^{4\pi|u^{*}|^2}<\infty.
$$
Moreover, using a standard argument we can prove that the functional $I$ is of class $C^1$ with
$$
I'(u)v=\int_{A_m}\nabla u\nabla v-\int_{A_m}f(u)v,~~~\forall u,v\in
H_0^1(A_m).
$$
Therefore, critical points of $I$ are precisely the weak solutions of $(P_m)$.

The next lemma ensures that the functional $I$ has the mountain pass geometry.

\begin{lemma}\label{mpg}

\begin{enumerate}

\item[$(a)$] There exist $r,\rho>0$ such that
$$
I(u)\geq\rho>0,~~ \mbox{for all}~~\|u\|_{A_m}=r.
$$
\item[$(b)$] There is $e\in H_0^1(A_m)$ such that
$$
\|e\|_{A_m}>r~~~\mbox{and}~~~ I(e)<0.
$$
\end{enumerate}

\end{lemma}

\noindent{\bf Proof.} Using the definition of $I$ and the growth of $f$, we obtain
$$
I(u)\geq\dfrac{1}{2}\int_{A_m}|\nabla
u|^2-\dfrac{\epsilon}{2}\int_{A_m}|u|^2-C\int_{A_m}|u|^qe^{\beta|u|^2},
$$
or equivalently,
$$
I(u)\geq\dfrac{1}{2}\int_{B}|\nabla
u^{*}|^2-\dfrac{\epsilon}{2}\int_{B}|u^{*}|^2-C\int_{B}|u^{*}|^qe^{\beta|u^{*}|^2}.
$$
By the Poincaré inequality,
$$
I(u)\geq\dfrac{1}{2}\int_{B}|\nabla
u^{*}|^2-\dfrac{\epsilon}{2\lambda_1}\int_{B}|\nabla
u^{*}|^2-C\int_{B}|u^{*}|^qe^{\beta|u^{*}|^2},
$$
where $\lambda_1$ is the first eigenvalue of $(-\Delta, H_0^1(B))$.
Fixing $\epsilon>0$ sufficiently small, we have
$C_1:=\dfrac{1}{2}-\dfrac{\epsilon}{2\lambda_1}>0$, from where it follows that
$$
I(u)\geq C_1\int_{B}|\nabla
u^{*}|^2-C\int_B|u^{*}|^qe^{\beta|u^{*}|^2}.
$$
Notice that, from Trudinger-Moser inequality (\ref{X1})
$$
e^{\beta|u^{*}|^2}\in L^2(B)
$$
and by continuous embedding
$$
|u^{*}|^q\in L^2(B).
$$
Since $H_0^1(B)\hookrightarrow L^{2q}(B)$ for all $q\geq 1$, by H\"older inequality
$$
\begin{array}{rl}
\displaystyle\int_B|u^{*}|^qe^{\beta|u^{*}|^2}&\leq\left(\displaystyle\int_B|u^{*}|^{2q}\right)^{1/2}\left(e^{2\beta|u^{*}|^2}\right)^{1/2}\\[0.3cm]
&\leq |u^{*}|_{2q,B}^q\left(\displaystyle\int_B e^{2\beta|u^{*}|^2}\right)^{1/2}\\[0.3cm]
&\leq C\|u^{*}\|_B^q\left(\displaystyle\int_B
e^{2\beta|u^{*}|^2}\right)^{1/2}.
\end{array}
$$

We claim that for $r>0$ small enough, we have
$$
\sup_{\|u^{*}\|_B=r}\int_Be^{2\beta|u^{*}|^2}<\infty.
$$
In fact, note that
$$
\int_Be^{2\beta|u^{*}|^2}=\int_Be^{2\beta\|u^{*}\|_B^2\left(\frac{|u^{*}|}{\|u^{*}\|_B}\right)^2}.
$$
Choosing $0<r\approx 0$ such that $\alpha:=2\beta r^2<4\pi$ and using the Trudinger-Moser inequality (\ref{X1}), 
$$
\sup_{\|u^{*}\|_B=r}\int_Be^{2\beta |u^{*}|^2}\leq\sup_{\|v\|_B\leq
1}\int_Be^{\alpha|v|^2}<\infty.
$$
Thus,
$$
I(u)\geq C_1\|u^{*}\|_B^2-C_2\|u^{*}\|_B^q.
$$
Fixing $q>2$, we derive
$$
I(u)\geq C_1r^2-C_2r^q:=\rho>0,
$$
for $r=\|u\|_{A_m}=\|u^{*}\|_B$ small enough, which shows that the item $(a)$ holds.

To prove $(b)$, first notice that

\noindent{\bf Claim 1.} For each $\epsilon>0$, there exists $\overline{s}_\epsilon>0$ such that
$$
F(s)\leq\epsilon f(s)s,~\mbox{for all}~ x\in A_m,~|s|\geq \overline{s}_\epsilon.
$$
In fact, from hypothesis $(H_1)$
$$
\left|\dfrac{F(s)}{sf(s)}\right|\leq\dfrac{M}{|s|},~\forall |s|\geq s_0.
$$
For $p>2$, the claim $1$ with $\epsilon=1/p>0$, guarantees the existence of $\overline{s}_\epsilon>0$ such that
$$
pF(s)\leq f(s)s,~\forall s\geq \overline{s}_\epsilon,
$$
which implies the existence of constant $C_1,C_2>0$ verifying
$$
F(s)\geq C_1|s|^{p}-C_2,~\forall s\geq 0.
$$
Thus, fixing $\varphi\in C_0^\infty(A_m)$ with $\varphi\geq 0$ and
$\varphi\neq 0$. For $t\geq 0$, we have
$$
\begin{array}{rl}
\displaystyle\int_{A_m}F(t\varphi)&\geq\displaystyle\int_{A_m}\left(C_1|t\varphi|^{p}-C_2\right)\\[0.3cm]
&\geq C_1|t|^{p}\displaystyle\int_{A_m}|\varphi|-C_2|A_m|,
\end{array}
$$
from where it follows that
\begin{equation}\label{4.1}
\int_{A_m}F(t\varphi)\geq C_3|t|^{p}-C_4.
\end{equation}
From $(\ref{4.1})$, if $t\geq 0$, 
$$
I(t\varphi)\leq
\dfrac{t^2}{2}\|\varphi\|_{A_m}^2-C_3|t|^{p}+C_4.
$$
Since $p>2$,
$$
I(t\varphi)\to-\infty,~~~\mbox{as}~~~t\to +\infty.
$$
Fixing $t_0\approx+\infty$ and let $e=t_0\varphi$, we get
$$
\|e\|_{A_m}\geq r~~~\mbox{and}~~~I(e)<0.~~\hbox{\fim}
$$

%%%%%%%%%%%%%%%%%%%%%%%%%%
%%%
%%% Condição (PS)
%%%
%%%%%%%%%%%%%%%%%%%%%%%%%%

The next lemma is crucial to prove that the energy functional $I$ satisfies the Palais-Smale condition and its proof can be found in \cite{DOR}.
\begin{lemma}\label{1.6}
Let $\Omega\subset\R^N$ be a bounded domain and $(u_n)$ be a sequence of functions in $L^1(\Omega)$ such that $u_n$ converging to $u\in L^1(\Omega)$ in $L^1(\Omega)$. Assume that $f(u_n(x))$ and $f(u(x))$ are also $L^1$ functions. If
$$
\int_\Omega|f(u_n)u_n|\leq C,~~~\mbox{para todo}~~~n\in\mathbb{N},
$$
then $f(u_n)$ converges in $L^1(\Omega)$ to $f(u)$.

\end{lemma}

\begin{lemma}\label{ld}
The functional $I$ satifies $(PS)_d$ condition, for all
$d\in(0,1/2)$.
\end{lemma}

\noindent {\bf Proof.} Let $d<1/2$ and $(u_n)$ be a $(PS)_d$ sequence for the functional $I$, i.e.,
$$
I(u_n)\to d~~~ \mbox{and}~~~I'(u_n)\to 0, ~\mbox{as}~ n\to+\infty.
$$
For each $n\in\mathbb{N}$, let us define $\epsilon_n=\displaystyle\sup_{\|v\|\geq 1}\{|I'(u_n)v|\}$, then
$$
|I'(u_n)v|\leq \epsilon_n\|v\|_m,
$$
for all $v\in H_0^1(A_m)$, where $\epsilon_n=o_n(1)$. Thus
\begin{equation}\label{d1}
\dfrac{1}{2}\int_{A_m}|\nabla u_n|^2-\int_{A_m}F(u_n)=d+o_n(1),~~~\forall n\in\mathbb{N}
\end{equation}
and
\begin{equation}\label{d2}
\left|\int_{A_m}\nabla u_n\nabla v-\int_{A_m}f(u_n)v\right|\leq\|v\|_m\epsilon_n,~~ \forall n\in\mathbb{N},~~v\in H_0^1(A_m).
\end{equation}
From $(\ref{d1})$ and Claim $1$, for any $\epsilon>0$, there is $n_0\in\mathbb{N}$ such that
$$
\dfrac{1}{2}\|u_n\|_m^2=\dfrac{1}{2}\int_{A_m}|\nabla u_n|^2\leq\epsilon+d+\int_{A_m}F(u_n)\leq C_\epsilon+\epsilon\int_{A_m}f(u_n)u_n,
$$
for all $n\geq n_0$. Using $(\ref{d2})$ with $v=u_n$, we get
$$
\left(\dfrac{1}{2}-\epsilon\right)\|u_n\|^2_m\leq C_\epsilon+\epsilon\|u_n\|_m,~\forall n\geq n_0.
$$
Thus, the sequence $(u_n)$ is bounded. Since $H_0^1(A_m)$ be a reflexive Banach space, there exits $u\in H_0^1(A_m)$ such that, for some subsequence,
$$
u_n\rightharpoonup u~~~\mbox{in}~~~ H_0^1(A_m).
$$
Furthermore, from compact embedding, 
$$
u_n\to u~~~\mbox{in}~~~ L^q(A_m),~q\geq 1
$$
and
$$
u_n(x)\to u(x)~~~\mbox{a.e. in}~~~ A_m.
$$
On the other hand, using $(\ref{d2})$ with $v=u_n$, we get
$$
-\epsilon_n\|u_n\|_m\leq\int_{A_m}|\nabla u_n|^2-\int_{A_m}f(u_n)u_n,
$$
which implies
$$
\int_{A_m}f(u_n)u_n\leq\|u_n\|_m^2-\epsilon_n\|u_n\|_m\leq C, ~~\forall n\in\mathbb{N}.
$$
From Lemma \ref{1.6}, $f(u_n)\to f(u)$ in $L^1(A_m)$. Then, there is $h\in L^1(A_m)$ such that
$$
|f(u_n(x))|\leq h(x),~\mbox{a.e. in}~ A_m,
$$
and from $(H_1)$, 
$$
|F(u_n)|\leq M h(x),~\mbox{a.e. in}~ A_m.
$$
Furthermore,
$$
F(u_n(x))\to F(u(x))~~\mbox{a.e. in}~~ A_m.
$$
Consequently, by the Lebesgue's dominated convergence,
$$
\int_{A_m}F(u_n)-\int_{A_m}F(u)=o_n(1).
$$
Thus, from $(\ref{d1})$,
$$
\dfrac{1}{2}\|u_n\|_m^2-\int_{A_m}F(u)-d=o_n(1),
$$
which implies,
\begin{equation}\label{d3}
\lim_{n\to\infty}\|u_n\|_m^2=2\left(d+\int_{A_m}F(u)\right).
\end{equation}
Using again $(\ref{d2})$ with $v=u_n$, we obtain
$$
\left|\|u_n\|_m^2-\int_{A_m}f(u_n)u_n\right|\leq o_n(1),
$$
from where we derive
$$
\begin{array}{lll}
\left|\displaystyle\int_{A_m}f(u_n)u_n-2\left(d+\displaystyle\int_{A_m}F(u)\right)\right|&\leq
\left|\|u_n\|_m^2-\displaystyle\int_{A_m}f(u_n)u_n\right|\\[0.3cm]
&+\left|\|u_n\|_m^2-2\left(d+\displaystyle\int_{A_m}F(u)\right)\right|.
\end{array}
$$
Then,
$$
\lim_{n\to\infty}\int_{A_m}f(u_n)u_n=2\left(d+\int_{A_m}F(u)\right).
$$
Furthermore, from $(H_2)$,
$$
\begin{array}{ll}
2\displaystyle\int_{A_m}F(u)&\leq 2\displaystyle\lim_{n\to\infty}\displaystyle\int_{A_m}F(u_n)\\[0.3cm]
&\leq\displaystyle\lim_{n\to\infty}\displaystyle\int_{A_m}f(u_n)u_n=2d+2\displaystyle\int_{A_m}F(u),
\end{array}
$$
which implies that $d\geq 0$.

\noindent{\bf Claim 2.} For any $v\in H_0^{1}(A_m)$, 
$$
\int_{A_m}\nabla u\nabla v=\int_{A_m}f(u)v.
$$
In fact, let us fix $v\in H_0^{1}(A_m)$ and notice that
$$
\begin{array}{ll}
\left|\displaystyle\int_{A_m}\!\!\!\nabla u\nabla v-\displaystyle\int_{A_m}\!\!\!f(u)v\right|\!\!\!&\leq\left|\displaystyle\int_{A_m}\!\!\!\nabla u_n\nabla v-\displaystyle\int_{A_m}\!\!\!\nabla u\nabla v\right|+\left|\displaystyle\int_{A_m}\!\!\!f(u_n)v-\displaystyle\int_{A_m}\!\!\!f(u)v\right|\\[0.3cm]
&+\left|\displaystyle\int_{A_m}\!\!\!\nabla u_n\nabla v-\displaystyle\int_{A_m}\!\!\!f(u_n)v\right|.
\end{array}
$$
Using Lemma \ref{1.6}, the weak convergence $u_n\rightharpoonup u$ in $H_0^1(A_m)$ and the estimate in $(\ref{d2})$, we derive
$$
\left|\int_{A_m}\nabla u\nabla v-\int_{A_m}f(u) v\right|\leq o_n(1)+\|v\|o_n(1),
$$
and the proof of Claim 2 is complete.

Notice that from $(H_2)$ and Claim 2,
$$
J(u)\geq\dfrac{1}{2}\int_{A_m}|\nabla u|^2-\dfrac{1}{2}\int_{A_m}f(u)u=0.
$$

Now, We split the proof into three cases:

\noindent{\bf Case 1.} The level $d=0$. By the lower semicontinuity of the norm, 
$$
\|u\|_m\leq\liminf_{n\to\infty}\|u_n\|_m,
$$
then
$$
\dfrac{1}{2}\|u\|_m^2\leq\dfrac{1}{2}\|u_n\|_m^2.
$$
Using $(\ref{d3})$,
$$
0\leq I(u)\leq\dfrac{1}{2}\liminf\|u_n\|^2-\int_{A_m}F(u),
$$
which implies
$$
0\leq I(u)\leq\int_{A_m}F(u)-\int_{A_m}F(u)=0,
$$
from where $I(u)=0$, or equivalently,
$$
\|u\|_m^2=2\int_{A_m}F(u).
$$
Using again $(\ref{d3})$, we derive
$$
\|u_n\|_m^2-\|u\|_m^2=o_n(1),
$$
since $H_0^1(A_m)$ be a Hilbert espace, 
$$
u_n\to u~~\mbox{in}~~ H_0^1(A_m).
$$
Therefore, $I$ verifies the Palais-Smale at the level $d=0$.

\noindent{\bf Case 2.} The level $d\neq 0$ and the weak limit $u\equiv 0$.

We will show that this can not occur for a Palais-Smale sequence.

\noindent{\bf Claim 3.} There are $q>1$ and a constant $C>0$ such that
$$
\int_{A_m}|f(u_n)|^q<C,~\forall n\in\mathbb{N}.
$$
In fact, from $(\ref{d3})$, for each $\epsilon>0$
$$
\|u_n\|_m^2\leq 2d+\epsilon, ~\forall n\geq n_0,
$$
for some $n_0\in\mathbb{N}$. Furthermore, from $(f_1)$,
$$
\int_{A_m}|f(u_n)|^q\leq C\int_{A_m}e^{4\pi qu_n^2}=C\int_{B}e^{4\pi\|u_n^{*}\|^2(\frac{u_n}{\|u_n^{*}\|})^2}.
$$
By the Trudinger-Moser inequality (\ref{X1}), the last integral in the equality above is bounded if $4\pi q\|u_n^{*}\|^2<4\pi$ and this occur if we take $q>1$ suficiently close to $1$ and $\epsilon$ small enough, because $d<1/2$, which proves the claim.

Then, using $(\ref{d2})$ with $v=u_n$, we obtain
$$
\left|\int_{A_m}|\nabla u_n|^2-\int_{A_m}f(u_n)u_n \right|\leq\epsilon_n\|u_n\|_m\leq \epsilon_nC,~ \forall n\in\mathbb{N}.
$$
Thus,
\begin{equation}\label{d4}
\|u_n\|_m^2\leq o_n(1)+\int_{A_m}f(u_n)u_n,~ \forall n\in\mathbb{N}.
\end{equation}
Furthermore, from H\"older inequality, we can estimate the integral above as follows
$$
\int_{A_m}f(u_n)u_n\leq\left(\int_{A_m}|f(u_n)|^q\right)^{1/q}\left(\int_{A_m}|u_n|^{q'}\right)^{1/q'},~\forall n\in\mathbb{N},
$$
and since $u_n\to 0$ in $L^{q'}(A_m)$, 
$$
\int_{A_m}f(u_n)u_n=o_n(1).
$$
Then, from $(\ref{d4})$,
\begin{equation}\label{vivo}
\|u_n\|_m^2\to 0,~\mbox{as}~n\to\infty,
\end{equation}
which contradicts $(\ref{d3})$, because
$$
\|u_n\|_m^2\to 2d\neq 0,~\mbox{as}~n\to\infty,
$$
proving that $d\neq 0$ and $u=0$ does not occur.

\noindent{\bf Case 3.} The level $d\neq 0$ and the weak limit $u\neq 0$. Since
$$
I(u)=\dfrac{1}{2}\|u\|_m^2-\int_{A_m}\!\!\!F(u)\leq\liminf_{n}\left(\dfrac{1}{2}\|u_n\|_m^2-\int_{A_m}\!\!\!F(u_n)\right)=d.
$$
we have $I(u)\leq d$.

\noindent{\bf Claim 4.} $I(u)=d$.

In fact, suppose by contradiction that $I(u)<d$, from definition of $I$, 
\begin{equation}\label{d5}
\|u\|_m^2<2\left(d+\int_{A_m}\!\!\!F(u)\right).
\end{equation}
On the other hand, if we consider the functions
$$
v_n=\dfrac{u_n^{*}}{\|u_n^{*}\|},~ n\in\mathbb{N}
$$
and
$$
v=u^{*}\left[2\left(d+\int_{B}\!\!F(u^{*})\right)\right]^{-1/2},
$$
we have  $\|v_n\|_B=1$ e $\|v\|_B<1$. Furthermore, since
$$
\int_B\!\!\!\nabla v_n\nabla \varphi=\|u_n\|^{-1}\!\!\!\int_{A_m}\!\!\!\nabla u_n\nabla\varphi\to\left[2\left(d+\!\!\int_B\!\!F(u^{*})\right)\right]^{-1/2}\!\!\!\int_B\!\!\nabla u\nabla \varphi=\int_B\!\!\!\nabla v\nabla\varphi,
$$
for every $\varphi\in C_0^\infty(B)$, i.e.,
$$
\int_B\nabla v_n\nabla\varphi-\int_B\nabla v\nabla\varphi=o_n(1),
$$
we have
$$
v_n\rightharpoonup v~~\mbox{in}~~H_0^1(B).
$$

\begin{claim}\label{viv35} There are $q>1$ and $n_0\in\mathbb{N}$ such that
$$
\int_{A_m}|f(u_n)|^q<C,~\forall n\geq n_0.
$$
\end{claim}

To prove this, we need the following result due to P.L. Lions \cite{L}.

\begin{proposition}\label{c2}
Let $(u_n)$ be a sequence in $H_0^1(\Omega)$ such that $|\nabla u_n|_2=1$ for all $n\in\mathbb{N}$. Furthermore, suppose that $u_n\rightharpoonup u$ in $H_0^1(\Omega)$ with $|\nabla u|_2<1$. If $u\neq 0$, then for each $1<p<\dfrac{1}{1-|\nabla u|_2^2}$, we have
$$
\sup_{n\in\mathbb{N}}\int_{\Omega}e^{4\pi pu_n^2}<\infty.
$$
\end{proposition}

From hypothesis $(f_1)$,
\begin{equation}\label{viv31}
\int_{A_m}|f(u_n)|^q\leq C\int_{A_m}e^{4\pi q u_n^2}=C\int_Be^{4\pi q\|u_n^{*}\|^2v_n^2}.
\end{equation}
The last integral in the above expression  is bounded. In fact, by Proposition \ref{c2}, it is suffices to prove that there are $q,~p>1$ and $n_0\in\mathbb{N}$ such that
\begin{equation}\label{viv30}
 q\|u_n^{*}\|^2\leq p<\dfrac{1}{1-\|v\|^2},~~\forall n\geq n_0.
\end{equation}
To prove that (\ref{viv30}) occur, notice that $I(u)\geq 0$ and $d<1/2$, which implies that
$$
2<\dfrac{1}{d-I(u)},
$$
from where it follows that
$$
2\left(d+\int_BF(u^{*})\right)<\dfrac{d+\displaystyle\int_BF(u^{*})}{d-I(u)}=\dfrac{1}{1-\|v\|_B^2}.
$$
Thus, for $q>1$ sufficiently close to $1$,
$$
2q\left(d+\int_BF(u)\right)<\dfrac{1}{1-\|v\|_B^2}.
$$
From $(\ref{d3})$, there are $p>1$ and $n_0\in\mathbb{N}$ such that
$$
q\|u_n^{*}\|^2\leq p<\dfrac{1}{1-\|v\|_B^2},~\forall n\geq n_0.
$$
Thus, (\ref{viv30}) occur. Therefore, Claim \ref{viv35} holds.

Now, we will show that $u_n\to u$ in $H_0^1(A_m)$. First, notice that from H\"older inequality and (\ref{viv35}),
$$
\int_{A_m}\!\!\!f(u_n)(u_n-u)\leq\int_{A_m}\!\!\!\left(|f(u_n)|^q\right)^{1/q}\left(\int_{A_m}|u_n-u|^{q'}\right)^{1/q'}\leq C|u_n-u|_{q',A_m},
$$
where $1/q+1/q'=1$. Since $u_n\to u$ in $L^{q'}(A_m)$, 
\begin{equation}\label{viv40}
\int_{A_m}\!\!\!f(u_n)(u_n-u)=o_n(1).
\end{equation}
Using $(\ref{d2})$ with $v=u_n-u$ and (\ref{viv40}), we obtain $\langle u_n-u,u_n\rangle=o_n(1)$, and since $u_n\rightharpoonup u$ in $H_0^1(A_m)$, 
$$
\|u_n-u\|_m^2=\langle u_n-u,u_n\rangle-\langle u_n-u,u\rangle=o_n(1).
$$
Then, $\|u_n\|_m^2\to\|u\|_m^2$ and this together with $(\ref{d3})$ contradicts $(\ref{d5})$. Which proves that $I(u)=d$, i.e.,
$$
\|u\|_m^2=2\left(d+\int_{A_m}\!\!\!F(u)\right).
$$
Furthermore, from $(\ref{d3})$, $\|u_n\|_m\to\|u\|_m$ as $n\to \infty$. Therefore
$$
u_n\to u~~\mbox{em}~~ H_0^1(A_m).~\hbox{\fim}
$$

From Lemma \ref{mpg} and the Mountain pass Theorem without compactness conditions (see \cite{W}), there is a $(PS)_{c_m}$ sequence $(u_n)\subset
H_0^1(A_m)$ such that
$$
I(u_n)\to c_m~~~\mbox{and}~~~I'(u_n)\to 0,
$$
where
$$
c_m=\inf_{\gamma\in\Gamma}\max_{t\in[0,1]} I(\gamma(t))
$$
and
$$
\Gamma=\{\gamma\in C([0,1],H_0^1(A_m)):~ \gamma(0)=0~ \mbox{and}~I(\gamma(1))<0\}.
$$
To conclude the proof of existence of positive solution for $(P)_m$, it remains to show that $c_m\in(-\infty,1/2)$. For this, we introduce the following Moser's functions (see \cite{M}):
$$
\overline{w}_n(x)=\dfrac{1}{\sqrt{2\pi}}\left\{
\begin{array}{cc}
\left(ln(n)\right)^{1/2},& 0\leq|x|\leq 1/n\\[0.3cm]
\dfrac{ln\dfrac{1}{|x|}}{\left(ln(n)\right)^{1/2}},& 1/n\leq|x|\leq 1\\[0.7cm]
0,&|x|\geq 1
\end{array}\right.
$$
Let $d_m>0$ and $x_m\in A_m$ such that $B_{d_m}(x_m)\subset A_m$ and define
$$
w_n(x)=\overline{w}_n\left(\dfrac{x-x_m}{d_m}\right),
$$
we have $w_n\in H_0^1(A_m)$, $\|w_n\|_{A_m}=1$ and $supp~ w_n\subset
B_{d_m}(x_m)$.

We claim that the exists $n\in\mathbb{N}$ such that
$$
\max_{t\geq 0}I(tw_n)<\dfrac{1}{2}.
$$
In fact, suppose by contradition that this is not the case. Then, there exist $t_n>0$ such that
\begin{equation}\label{e1}
\max_{t\geq 0}I(tw_n)=I(t_nw_n)\geq \dfrac{1}{2}.
\end{equation}
It follows from $(\ref{e1})$ and $(H_1)$ that
\begin{equation}\label{e2}
t_n^2\geq 1.
\end{equation}
Furthermore, $\dfrac{d}{dt}I(tw_n)\left|_{t=t_n}\right.=0$, i.e.,
\begin{equation}\label{e5}
t_n^2=\int_{A_m}f(t_nw_n)t_nw_n,
\end{equation}
which implies that
\begin{equation}\label{e3}
t_n^2\geq\int_{B_{d_m/n}(x_m)}f(t_nw_n)t_nw_n.
\end{equation}
In what follows, we fix a positive constant $\beta_m$ verifying
\begin{equation}\label{eq7}
\beta_m>\dfrac{1}{2\pi d_m^2}.
\end{equation}
From $(H_3)$, there exists $s_m=s_m(\beta_m)>0$ such that
\begin{equation}\label{bm}
f(s)s\geq\beta_m e^{4\pi s^2},~\forall s\geq s_m.
\end{equation}
Using $(\ref{bm})$ in $(\ref{e3})$ and the definition of $w_n$ in $B_{d_m/n}(0)$, we obtain
\begin{equation}
t_n^2\geq \beta_m\pi \dfrac{d_m^2}{n^2}e^{2t_n^2ln(n)}
\end{equation}
for $n$ large enough, or equivalently,
\begin{equation}\label{1}
t_n^2\geq \beta_m\pi d_m^2e^{2ln(n)(t_n^2-1)},
\end{equation}
it implies that the sequence $(t_n)$ is bounded. Moreover, from $(\ref{1})$ and
$(\ref{e2})$,
$$
t_n^2\to 1,~\mbox{as}~ n\to \infty.
$$
Now, let us define
$$
C_n=\{x\in B_{d_m}(x_m):t_nw_n(x)\geq s_m\}
$$
and
$$
D_n=B_{d_m}(x_m)\setminus C_n.
$$
With the above notations and using $(\ref{e5})$,
$$
t_n^2\geq\int_{B_{d_m/n}(x_m)}f(t_nw_n)t_nw_n=\int_{C_n}f(t_nw_n)t_nw_n+\int_{D_n}f(t_nw_n)t_nw_n
$$
and by $(\ref{bm})$,
$$
t_n^2\geq\int_{D_n}f(t_nw_n)t_nw_n+\beta_m\int_{C_n}e^{4\pi
t_n^2w_n^2}
$$
or equivalently,
\begin{equation}\label{e2.2}
t_n^2\geq\int_{D_n}f(t_nw_n)t_nw_n+\beta_m\int_{B_{d_m}(x_m)}e^{4\pi
t_n^2w_n^2}-\beta_m\int_{D_n}e^{4\pi t_n^2w_n^2}.
\end{equation}

Notice that
$$
w_n(x)\to 0 ~~\mbox{a.e. in}~~ B_{d_m}(x_m),
$$
$$
\chi_{D_n}(x)\to 1~~\mbox{a.e. in}~~ B_{d_m}(x_m)
$$
and
$$
e^{4\pi t_n^2w_n^2}\chi_{D_n}\leq e^{4\pi t_n^2s_m^2}\in L^1(B_{d_m}(x_m)).
$$
Then, by Lebesgue's dominated convergence
\begin{equation}\label{e2.1}
\lim_n\int_{D_n}e^{4\pi t_n^2w_n^2}=\lim_n\int_{B_{d_m}(x_m)}e^{4\pi
t_n^2w_n^2}\chi_{D_n}=\int_{B_{d_m}(x_m)}1=\pi d_m^2.
\end{equation}
Furthermore,
$$
f(t_nw_n)t_nw_n\chi_{D_n}\leq Ct_nw_ne^{4\pi t_n^2w_n^2}\leq Cs_m
e^{4\pi s_m^2}\in L^1(B_{d_m}(x_m))
$$
and
$$
f(t_nw_n(x))t_nw_n(x)\chi_{D_n}(x)\to 0~~\mbox{a.e. in}~~B_{d_m}(x_m).
$$
Thus, using again Lebesgue's dominated convergence,
\begin{equation}\label{e1.1}
\lim_n\int_{D_n}f(t_nw_n)t_nw_n=0
\end{equation}
Passing to the limit $n\to\infty$ in $(\ref{e2.2})$ and using
$(\ref{e2.1})$ and $(\ref{e1.1})$,
$$
1\geq\beta_m\lim_n\int_{B_{d_m}(x_m)}e^{4\pi t_n^2w_n^2}-\beta_m\pi
d_m^2.
$$
Since $t_n^2\geq 1$, we get
\begin{equation}\label{e2.3}
1\geq\beta_m\lim_n\left[\int_{B_{d_m}(x_m)}e^{4\pi w_n^2}\right]-\beta_m\pi
d_m^2.
\end{equation}
On the other hand, since
$$
\int_{B_{d_m}(x_m)}e^{4\pi
w_n^2}=d_m^2\int_{B_1(0)}e^{4\pi\overline{w}_n^2}=d_m^2\left\{\dfrac{\pi}{n^2}e^{4\pi\frac{1}{2\pi}ln(n)}+2\pi\int_{1/n}^1e^{4\pi\frac{1}{2\pi}\frac{[ln(1/r)]^2}{ln(n)}}rdr\right\},
$$
making a changing of variables $s=ln(1/r)/ln(n)$,
$$
\int_{B_{d_m}(x_m)}e^{4\pi
w_n^2}=\pi d_m^2+2\pi
d_m^2ln(n)\int_0^1e^{2s^2ln(n)-2sln(n)},
$$
and since
$$
\lim_{n\to\infty}\left[ 2ln(n)\int_0^1e^{2ln(n)(s^2-s)}ds\right]=2,
$$
we have
$$
\lim_{n\to\infty}\int_{B_{d_m}(x_m)}e^{4\pi w_n^2}=\pi d_m^2+2\pi d_m^2=3\pi d_m^2.
$$
Using the last limit in $(\ref{e2.3})$, we get
$$
1\geq 3\beta_m\pi d_m^2-\beta\pi d_m^2=2\beta\pi d_m^2,
$$
from where we derive
$$
\beta_m\leq\dfrac{1}{2\pi d_m^2},
$$
which contradicts the choice of $\beta_m$ in $(\ref{eq7})$. Then,
$$
\max_{t\geq 0}I(tw_n)<\dfrac{1}{2},
$$
proving that $c_m<1/2$, for any $m\in\mathbb{N}$ fixed arbitrarily. \fim

\section{Proof of Theorem \ref{sb}.}

To  proof Theorem\ref{sb}, we will use the following proposition.
\begin{proposition}\label{p4u}
Let $A$ be an angular sector contained on the positive semiplane of $\R^2$ such that one of its boundary lies in $x_2$ axis, and denote such boundary of $A$ by $B_0=\{x=(x_1,x_2)\in A:~ x_2=0\}$. Consider $A'$ the reflection $A$ with respect to $x_2$ axis. Suppose that $u$ is a solution of the following problem:
$$
\left\{
\begin{array}{cl}
-\Delta u=f(u),& \mbox{ in }~~ A, \\[0.2cm]
u=0,&\mbox{on}~~ B_0,
\end{array}
\right.\leqno{(P)}
$$
where $f$ is a real, continuous and odd function. Then, the function $\tilde{u}$ such that $\tilde{u}= u$ in $A$  and $\tilde{u}$ is antisymmetric
with respect to $x_2$ axis,
$$
\tilde{u}(x_1,x_2)=\left\{
\begin{array}{ll}
u(x_1,x_2),&\mbox{in}~~A\\[0.2cm]
-u(x_1,-x_2),&\mbox{in}~~A'\\[0.2cm]
0,&\mbox{on}~~B_0
\end{array}
\right.
$$
satisfies
$$
-\Delta \tilde{u}=f(\tilde{u})~~~\mbox{in}~~~A\cup A'.
$$
\end{proposition}

\noindent{\bf Proof.} Since $u$ be a solution of $(P)$, we have
$$
\int_{A}\nabla u\nabla
\varphi=\int_{A}f(u)\varphi,~~\forall \varphi\in
C^\infty_c(A).
$$
we want to prove that
$$
\int_{A\cup A'}\nabla \tilde{u}\nabla
\phi=\int_{A\cup A'}f(\tilde{u})\phi,~\forall \phi\in
C^\infty_0(A\cup A').
$$
For any $\phi\in C^\infty_0(A\cup A')$,
$$
\int_{A\cup A'}f(\tilde{u})\phi=\int_{A}f(u(x_1,x_2))\phi(x_1,x_2)+\int_{A'}f(-u(x_1,-x_2))\phi(x_1,x_2).
$$
Since $f$ be an odd function, 
$$
\begin{array}{ll} \displaystyle\int_{A\cup A'}f(\tilde{u})\phi
&=\displaystyle\int_{A}f(u(x_1,x_2))\phi(x_1,x_2)+\int_{A'}f(-u(x_1,-x_2))\phi(x_1,x_2)\\[0.3cm]
&=\displaystyle\int_{A}f(u(x_1,x_2))\phi(x_1,x_2)-\int_{A'}f(u(x_1,-x_2))\phi(x_1,x_2)\\[0.3cm]
&=\displaystyle\int_{A}f(u(x_1,x_2))\phi(x_1,x_2)-\int_{A}f(u(x_1,x_2))\phi(x_1,-x_2).
\end{array}
$$
Thus
\begin{equation}\label{e11}
\int_{A\cup A'}f(\tilde{u})\phi=\int_{A}f(u)\psi,
\end{equation}
where $\psi(x_1,x_2)=\phi(x_1,x_2)-\phi(x_1,-x_2)$. On the other hand,
$$
\begin{array}{ll}
\displaystyle\int_{A\cup A'}\nabla\tilde{u}\nabla\phi
&=\displaystyle\int_{A}\nabla u(x_1,x_2)\nabla\phi(x_1,x_2)-\int_{A'}\nabla u(x_1,-x_2)\nabla\phi(x_1,x_2)\\[0.3cm]
&=\displaystyle\int_{A}\nabla u(x_1,x_2)\nabla\phi(x_1,x_2)-\int_{A}\nabla u(x_1,x_2)\nabla (\phi(x_1,-x_2))\\[0.3cm]
&=\displaystyle\int_{A}\nabla u(x_1,x_2)\nabla(\phi(x_1,x_2)-\phi(x_1,-x_2)).
\end{array}
$$
Then,
\begin{equation}\label{e12}
\int_{A\cup A' }\nabla\tilde{u}\nabla\phi=\int_{A}\nabla
u\nabla\psi.
\end{equation}
The function $\psi$ does not in general belong to $C^\infty_0(A)$. Therefore, $\psi$ can not be used as a function
test (in the definition of weak solution on $H^1(A)$). On the other hand, if we consider the sequence of functions $(\eta_k)$ in $C^{\infty}(\R)$, definided by
$$
\eta_k(t)=\eta(kt),~t\in\R,~k\in\mathbb{N},
$$
where $\eta\in C^\infty(\R)$ is a function such that
$$
\eta(t)=\left\{
\begin{array}{ll}
0, & \mbox{if}~~t<1/2,\\
1, & \mbox{if}~~t>1.
\end{array}
\right.
$$
Then
$$
\varphi_k(x_1,x_2):=\eta_k(x_2)\psi(x_1,x_2)\in C^\infty_0(A),
$$
which implies that
\begin{equation}\label{e14}
\int_{A}\nabla
u\nabla\varphi_k=\int_{A}f(u)\varphi_k,~k\in\mathbb{N}.
\end{equation}
From (\ref{e11}), (\ref{e12}) and (\ref{e14}), we can conclude de proof, in view of the following limits
$$
(I)~~~~~~~~\int_{A}\nabla u\nabla\varphi_k\rightarrow
\int_{A}\nabla u\nabla\psi
$$
e
$$
(II)~~~~~~~~\int_{A}f(u)\varphi_k\rightarrow \int_{A}f(u)\psi,
$$
as $k\rightarrow \infty$ occur. To see that $(I)$ occur, notice that
$$
\int_{A}\nabla u\nabla\varphi_k=\int_{A}\eta_k\nabla u\nabla \psi+\int_{A}\dfrac{\partial u}{\partial x_2}k\eta'(kx_2)\psi.
$$
Clearly,
$$
\int_{A}\eta_k\nabla u\nabla \psi\to\int_{A}\nabla u\nabla \psi,~\mbox{as}~~k\to\infty.
$$
Assim, resta mostrar que
\begin{equation}\label{stR}
\int_{A}\dfrac{\partial u}{\partial x_2}k\eta'(kx_2)\psi\to 0~~\mbox{as}~~k\to\infty.
\end{equation}
In fact,
$$
\left|\int_{A}\dfrac{\partial u}{\partial x_2}k\eta'(kx_2)\psi\right|\leq kMC\int_{0<x_2<1/k}\left|\dfrac{\partial u}{\partial x_2}\right|x_2\leq MC\int_{0<x_2<1/k}\left|\dfrac{\partial u}{\partial x_2}\right|,
$$
where $C=\displaystyle\sup_{t\in[0,1]}|\eta'(t)|$ and $M>0$ is such that
$$
|\psi(x_1,x_2)|\leq M|x_2|,~ \forall (x_1,x_2)\in A\cup A',
$$
and since
$$
\int_{0<x_2<1/k}\left|\dfrac{\partial u}{\partial x_2}\right|\to 0,~~\mbox{as}~~k\to\infty,
$$
the limit in (\ref{stR}) occur. The item $(II)$ is an imediatly consequence of the Lebesgue's dominated convergence.

Now, for each $m\in\mathbb{N}$, we apply the Proposition \ref{p4u} to the solution $u$ of problem $(P)_m$. Let $A_m'$ be the reflection of $A_m$ in one of its sides. On $A_m\cup A_m'$, wecan the define the function $\tilde{u}$ such that $\tilde{u}=u$ on $A_m$, and $\tilde{u}$ is antisymmetric with respect to the side of reflection. Now, let $A_m''$ be the reflection of $A_m\cup A_m'$ in one of its sides and
$\tilde{\tilde{u}}$ the function defined on $A_m\cup A_m'\cup A_m''$
such that $\tilde{\tilde{u}}=\tilde{u}$ on $A_m\cup A_m'$ and
$\tilde{\tilde{u}}$ is antisymmetric with respect to the side of reflection.Repetindo este procedimento, após um número finito de reflexões. Repeating this procedure, after finite steps, we finally obtain a function defined on the whole unit ball $B$,  denoted by $u_m$. Clearly, $u_m$ satisfies the Dirichlet condition on the boundary $\partial B$. That is, $u_m$ is a sign-changing solution of problem $(P)$. Since for every $m\in\mathbb{N}$, problem $(P)_m$ admits a positive solution, hence there exist infinitely many differents sign-changing solutions, and the proof of Theorem  \ref{sb} is complete. \fim

\vspace{0.5cm}

\begin{figure}[!htb]
\begin{center}
\includegraphics[scale=0.25]{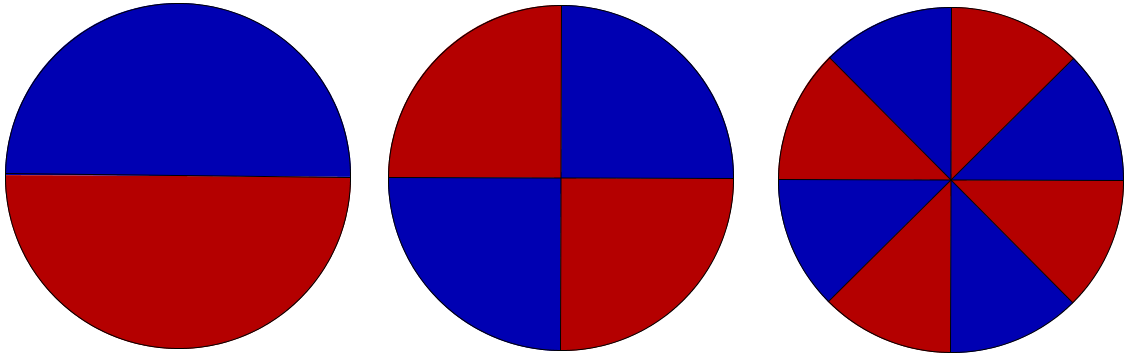}
\end{center}
\caption{Signal of solutions}
\label{ddd}
\end{figure}

\noindent {\bf Remark.} In figure $\ref{ddd}$, we represent the signal of three solutions, corresponding to the cases $m=1,~ m=2~\mbox{and}~ m=3$, respectively. The blue color represents the regions where the solutions are negative and the red color, the regions where the solutions are positive.

We show in Figure $3$ the profile of solution for the case $m=2$.

\begin{figure}[!htb]
\begin{center}
\includegraphics[scale=0.8]{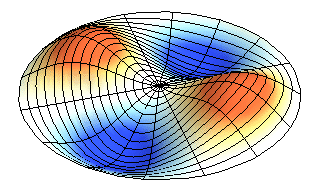}
\end{center}
\caption{Case $m=2$}
\end{figure}

\begin{remark} It is possible to make a version of Theorem \ref{sb} with Neumann boundary condition using the same arguments that we used here, but we have to work with another version of Trudinger-Moser inequality in $H^1(\Omega)$ due to Adimurthi-Yadava \cite{AYa}, which says that if $\Omega$ is a bounded domain with smooth boundary, then for any $u\in H^1(\Omega)$,
\begin{equation} \label{X02}
\int_{\Omega}
e^{\alpha u^2}< +\infty, \,\,\,\, \mbox{ for all }\,\,\alpha >0.
\end{equation}
Furthermore, there exists a positive constant $C=C(\alpha,|\Omega|)$ such that
\begin{equation} \label{X11}
\sup_{||u||_{H^{1}(\Omega)} \leq 1} \int_{\Omega} e^{\alpha u^2} \leq C , \,\,\,\,\,\,\, \forall \, \alpha  \leq 2 \pi .
\end{equation} 

\end{remark}
\vspace{0.3 cm}

\noindent \textbf{Acknowledgments.}  The author would like to thank professor Claudianor O. Alves for his excellent advices during the graduate school of him.

\end{document}